\documentclass[reqno,11pt]{amsart}

\usepackage{amsmath}
\usepackage{amssymb}
\usepackage{amsthm}
\usepackage{fullpage}
\usepackage{bm}
\usepackage{hyperref}
\usepackage{graphicx}
\usepackage{enumitem}

\usepackage[toc]{appendix}
  
\usepackage{pgfplots}
\pgfplotsset{compat=newest}
\usetikzlibrary{plotmarks}
\usetikzlibrary{arrows.meta}
\usetikzlibrary{calc}
\usepgfplotslibrary{patchplots}
\usepackage[font=small,labelfont=bf]{caption}

\usepackage{listings}
\usepackage{colortbl}
\newtheorem{thm}{Theorem}[section]
\newtheorem{lemma}[thm]{Lemma}
\newtheorem{prop}[thm]{Proposition}
\theoremstyle{remark}
\newtheorem{rmk}[thm]{Remark}

\def\RR{\mathbb{R}}

\begin{document}

\renewcommand{\thefootnote}{$\star$} 

\title[{Any three eigenvalues do not determine a triangle}]{Any three eigenvalues do not determine a triangle}

\author{Javier G\'omez-Serrano$^{*}$} 
\address[Corresponding author]{Department of Mathematics, Princeton University, Princeton, NJ 08544, USA}
\email{jg27@math.princeton.edu}

\author{Gerard Orriols}
\address{Cambridge University}
\email{go262@cam.ac.uk}

\thanks{}

\begin{abstract}
Despite the moduli space of triangles being three dimensional, we prove the existence of two triangles which are not isometric to each other for which the first, second and fourth Dirichlet eigenvalues coincide, establishing a numerical observation from Antunes--Freitas \cite{Antunes-Freitas:inverse-spectral-problem-triangles}. The two triangles are far from any known, explicit cases. To do so, we develop new tools to rigorously enclose eigenvalues to a very high precision, as well as their position in the spectrum. This result is also mentioned as (the negative) part of \cite[Conjecture 6.46]{Laugesen-Siudeja:triangles-survey}, \cite[Open Problem 1]{Grieser-Maronna:isospectral-triangle} and \cite[Conjecture 3]{Lu-Rowlett:spectral-amm}.
\end{abstract}

\maketitle

\setcounter{footnote}{0}
\renewcommand*{\thefootnote}{\arabic{footnote}}
\section{Introduction}


Mark Kac coined the term ``hearing the shape of a drum'' \cite{Kac:can-one-hear-shape-drum} in 1965 for the problem of the determination of a domain $D$ given the spectrum of the Laplacian. Since then, the spectral determination of (mostly planar) domains has become a fundamental question in geometric analysis. Throughout this paper we will work with Dirichlet boundary conditions, that is, we will consider the set of real numbers $0 < \lambda_1 < \lambda_2 \leq \lambda_3 \leq \ldots$ that solve
\begin{align}
 -\Delta u_k & = \lambda_k u_k  \text{ in } D \nonumber \\
u_k & = 0  \text{ on } \partial D. \label{laplacian-pde}
\end{align}

Most of the results in the literature are negative. In particular, for euclidean polygons, the first example of a pair of non-isometric polygons with the same spectrum is due to Gordon, Webb and Wolpert \cite{Gordon-Webb-Wolpert:isospectral-domains-surfaces-riemannian-orbifolds}. In the case of Riemannian manifolds, even the local geometry of isospectral manifolds can be different \cite{Szabo:isospectral-metrics}.

We now mention some of the very few positive results for the inverse spectral problem. For negatively curved manifolds, there are spectral rigidity \cite{Croke-Sharafutdinov:spectral-rigidity-negatively-curved-manifold,Guillemin-Kazhdan:inverse-spectral-negatively-curved-2-manifolds} and compactness results \cite{Osgood-Phillips-Sarnak:compact-isospectral-surfaces}, some of which have been recently extended to Anosov surfaces \cite{Paternain-Salo-Uhlmann:spectral-rigidity-anosov-surfaces}. See also \cite{Osgood-Phillips-Sarnak:compact-isospectral-plane-domains,Osgood-Phillips-Sarnak:compact-isospectral-plane-domains-annals} for a proof of these compactness results in the case of planar domains. Zelditch \cite{Zelditch:inverse-spectral-problem-analytic-annals} showed, by computing wave invariants, that any analytic bounded planar domain that has an even symmetry with respect to a line is spectrally determined. This was later generalized to higher dimensional domains \cite{Hezari-Zelditch:inverse-spectral-problem-Rn}. Without any symmetry and analyticity assumptions, only few results exist on the plane: a family of perturbed disks \cite{Watanabe:plane-domains-spectrally-determined}, semi-regular polygons in the class of convex piecewise smooth domains \cite{Enciso-GomezSerrano:spectral-semiregular} and very recently, Hezari and Zelditch \cite{Hezari-Zelditch:inverse-spectral-problem-ellipses} have shown that ellipses of small eccentricity are spectrally determined in the class of smooth domains, even without assumptions on convexity or closeness to the ellipse.

Instead of focusing on a wide class of domains, one could try to solve the inverse spectral problem in a smaller class. Simple, positive examples include the regular $n$-gon in the class of $n$-gons (due to the isoperimetric inequality) or rectangles in the class of rectangles (since the eigenvalues are explicit). A nontrivial result is to show that triangles are determined in the class of triangles (see \cite{Durso:phd-thesis} for a proof using the wave trace and a simpler one \cite{Grieser-Maronna:isospectral-triangle} that only uses the heat kernel and elementary calculations). Recently, trapezoids have been shown to be determined in the class of trapezoids (under Neumann boundary conditions) \cite{Hezari-Lu-Rowlett:neumann-isospectral-trapezoids}, as well as parallelograms and acute trapezoids \cite{Lu-Rowlett:spectral-amm} in their respective classes.

However, most known results use the full spectrum (i.e.~an infinite amount of quantities); the only non-explicit results using only finitely many eigenvalues rely on isoperimetric inequalities, most remarkably the Payne--P\'olya--Weinberger conjecture proved by Ashbaugh and Benguria~\cite{Ashbaugh-Benguria:ppw-conjecture}, which implies that the first two eigenvalues determine the disk among all domains of $\mathbb{R}^n$. See~\cite{Henrot:book-spectral-geometry} for a discussion on extremizers of functions of eigenvalues. Being in a finite dimensional ambient space (such as the space of $n$-gons for a fixed $n$) it is expected that a finite amount of eigenvalues should suffice to characterize an $n$-gon. Chang and DeTurck \cite{Chang-DeTurck:n-eigenvalues-triangle} showed that given a triangle $T$, there is a finite number $N(T)$ such that the first $N(T)$ eigenvalues distinguish $T$ from any other triangle. However, $N(T)$ is not known to be uniformly bounded. Since the moduli space of triangles is 3 dimensional, it is conceivable (see \cite{Laugesen-Siudeja:maximizing-neumann-triangles}) that 3 eigenvalues could characterize every triangle. Our paper is a first step in this direction by showing (on the negative) that if that is the case, not any 3 eigenvalues suffice. Specifically we prove the following theorem:

\begin{thm}
  \label{main_thm}
  There exist two triangles $T_\mathrm{A}$ and $T_\mathrm{B}$, not isometric to each other, such that $\lambda_i(T_\mathrm{A}) = \lambda_i(T_\mathrm{B})$, for $i=1,2,4$.
\end{thm}

Antunes and Freitas \cite{Antunes-Freitas:inverse-spectral-problem-triangles} had observed numerically the presence of a saddle point for $\lambda_4 / \lambda_1$ around which $\lambda_2 / \lambda_1$ is regular, which would imply the existence of such two triangles (see Figure~\ref{tworegions}). Moreover they conjectured that the three first eigenvalues $\lambda_1, \lambda_2$ and $\lambda_3$ do determine the shape of a triangle. This question has recently attracted a lot of attention. Laugesen and Siudeja mention both (positive and negative) directions in \cite[Conjecture 6.46]{Laugesen-Siudeja:triangles-survey}, as well as Grieser and Maronna \cite[Open Problem 1]{Grieser-Maronna:isospectral-triangle} and Lu and Rowlett \cite[Conjecture 3]{Lu-Rowlett:spectral-amm}. Despite being the simplest polygonal class, the spectral theory of triangles is far from well understood and many advances have been made recently, for example \cite{Judge-Mondal:hot-spots-triangles}, \cite{Judge-Hillariet:spectral-simplicity} and \cite{Ourmieres:flat-triangles}.

\begin{rmk}
  \label{rmk_open}
  With no extra effort, we also show that there exists an open set in the moduli space of triangles such that for every triangle from this set there is another one with the same eigenvalues $\lambda_i$, $i = 1, 2, 4$.
\end{rmk}


We do not know how to show the existence of a saddle point for the quotient of two eigenvalues in a region far away from the explicit triangles; that would require controlling up to two derivatives of the eigenvalues, which have very difficult expressions. Instead, we use a topological argument relying only on ``$C^0$ computations'' and using very tight, rigorous enclosures of eigenvalues, as well as a control of their position in the spectrum. While upper bounds are easy to obtain using Rayleigh--Ritz quotients over suitable test functions, guaranteed lower bounds (and indices) are much more difficult. 

There are many classical results in the literature \cite{Kato:upper-lower-bounds-eigenvalues,Lehmann:optimale-eigenwerte,Behnke-Goerisch:inclusions-eigenvalues,Fox-Henrici-Moler:approximations-bounds-eigenvalues,Moler-Payne:bounds-eigenvalues,Still:computable-bounds-eigenvalues} showing the existence of an eigenvalue close to an approximate one. The main philosophy is that if $(\lambda_\text{app},u_\text{app})$ is an eigenpair that satisfies the equation up to an error bounded by $\delta$, then there is a true eigenpair $(\lambda,u)$ at a distance $C\delta^k$. Recently, Barnett and Hassell \cite{Barnett-Hassell:quasi-ortogonality-dirichlet-eigenvalues} improved on the classical estimates by a factor of $O(\sqrt{\lambda})$ using quasi-orthogonality arguments. Unfortunately, none of these results can assert the position in the spectrum of the eigenvalue without knowing any a priori bounds.

In order to circumvent this issue and obtain the position, Plum \cite{Plum:eigenvalues-homotopy-method} proposed a homotopy method linking the eigenvalues of the domain of the problem with another, known domain (the base problem). See also the intermediate method \cite{Weinstein-Stenger:eigenvalues-book,Goerisch:stufenverfahren-eigenwerten,Beattie-Goerisch:lower-bounds-eigenvalues} for another example of connecting the problem to a known domain. The difficulty in our case is that we have to solve many eigenvalue problems and all of them are far from any known domain, yielding impractical times at the precision needed. In \cite{Liu-Oishi:verified-eigenvalues-laplacian-polygons} and \cite{Liu:framework-verified-eigenvalues}, Liu and Oishi found rigorous lower bounds in terms of the solutions of a finite dimensional problem given by a Finite Element Method discretization (see also \cite{Carstensen-Gedicke:lower-bounds-eigenvalues} for similar bounds). Again, at the level of the required precision for Theorem \ref{main_thm} the number of elements that we need in the mesh is too high.

In this paper we propose a combination of the two families of methods in two passes. In a first pass, we separate the first 4 eigenvalues from the rest of the spectrum (using the method of \cite{Liu:framework-verified-eigenvalues}). At that point the enclosures that we find are big and not admissible. In a second pass we find a very good approximation of 4 eigenvalues and eigenvectors below the threshold (which have to be necessarily close to $\lambda_1,\lambda_2,\lambda_3,\lambda_4$). Combining now the stability methods from Barnett--Hassell \cite{Barnett-Hassell:quasi-ortogonality-dirichlet-eigenvalues} with the a priori knowledge of the lower bound of $\lambda_5$ this yields very small enclosures of the first 4 eigenvalues and their respective positions in the spectrum. This new method combines the different strengths of the two families and their characters: the global problem (finding the order) versus the local problem (the refinement of its value).

In order to find a very accurate representation of the eigenvalue and eigenvector we will use the Method of Particular Solutions (MPS). This method was introduced by Fox, Henrici and Moler \cite{Fox-Henrici-Moler:approximations-bounds-eigenvalues} and has been later adapted by many authors (see \cite{Antunes-Valtchev:mfs-corners-cracks,Read-Sneddon-Bode:series-method-mps,Golberg-Chen:mfs-survey,Fairweather-Karageorghis:mfs-survey,Betcke:generalized-svd-mps} as a sample, and the thorough review \cite{Betcke-Trefethen:method-particular-solutions}). The main idea is to consider a set of functions that solve the eigenvalue problem without boundary conditions as a basis, and writing the solution of the problem with boundary conditions as a linear combination of them, solving for the coefficients that minimize the error on the boundary. Typically, the choices have been products of Bessel functions and trigonometric polynomials centered at certain points. The different choices for these points have a big impact in the performance of the method. Very recently, Gopal and Trefethen \cite{Gopal-Trefethen:new-laplace-solver-pnas,Gopal-Trefethen:laplace-solver-detailed} have developed a new way of selecting the base functions in such a way to yield root exponential convergence (the \textit{lightning Laplace solver}). We use their algorithm (though our own implementation in C++) to obtain accurate, fast approximations of the eigenpairs. We stress that these methods produce accurate approximations but there is no explicit control of the error with respect to the true solution. In order to overcome this issue and make all the bounds (such as the defect) rigorous we will use interval arithmetics to bound errors whenever needed.

By working this way, the method outlined above produces a rigorous enclosure of the eigenvalues corresponding to a particular triangle, but using only this information it is impossible to show that two given triangles share the same eigenvalues since we only obtain (narrow) bounds and we are looking at a closed condition. Instead, due to the stability of the problem, we can transform the closed condition into an open condition using a topological argument and check it by computing multiple bounds along the sides of a parallelogram in the moduli space of triangles, at the price of an increased computational cost. To mitigate this cost, we reduce the amount of points to be computed using stability estimates (with explicit constants) for the eigenvalues of a small open ball of triangles around a given one. We remark that perturbation methods will not work directly since the two triangles $T_\mathrm{A}$ and $T_\mathrm{B}$ are far from the triangles for which the spectrum is explicitly known.

In the recent years, the application of calculations done by computers to mathematical proofs has become more popular due to the increment of computational resources, but in order to make sure that their results are rigorous, we need to control the errors that floating point arithmetic can accumulate. This is usually done by means of interval arithmetic, in which the data that a computer stores for a real number is an interval (two endpoints, or a midpoint and a radius) of real numbers, stored by two high-precision floating point numbers, instead of just one.

Operations between intervals are implemented to return intervals which are guaranteed to contain every possible result when the operands belong to the input intervals. For example, if $[x] = [\underline{x}, \overline{x}]$ and $[y] = [\underline{y}, \overline{y}]$ are two intervals, their sum will can be given by the interval $[x] + [y] = [\underline{x} + \underline{y}, \overline{x} + \overline{y}]$ and their product by $[x] \cdot [y] = [\min\{\underline{x}\underline{y}, \underline{x}\overline{y}, \overline{x}\underline{y}, \overline{x}\overline{y}\}, \max\{\underline{x}\underline{y}, \underline{x}\overline{y}, \overline{x}\underline{y}, \overline{x}\overline{y}\}]$. The same rule applies to function implementations: a function $f$ evaluated on $[x]$ should return an interval containing every $f(x)$ for $x \in [x]$. We refer to the book \cite{Tucker:validated-numerics-book} for an introduction to validated numerics, in which most of the techniques used here are explained, and to the survey \cite{GomezSerrano:survey-cap-in-pde} and the recent book \cite{Nakao-Plum-Watanabe:cap-for-pde-book} for a more specific treatment of computer-assisted proofs in PDE.

The paper is organized as follows. In Section \ref{section:onetriangle} we explain how to get tight, rigorous enclosures of the eigenvalues corresponding to a single triangle. Section \ref{section:twotriangles} is devoted to the extension of the previous method to a continuous set of triangles, and to the proof of Theorem \ref{main_thm}. Finally, in Section \ref{section:implementation} we give some implementation details.

\section{Enclosing the spectrum of one triangle}
\label{section:onetriangle}
In this section we explain how we obtain rigorous bounds for the eigenvalue quotients $\xi_{21} = \frac{\lambda_2}{\lambda_1}, \xi_{41} = \frac{\lambda_4}{\lambda_1}$ for a particular triangle, using the two passes explained above. This will be the main ingredient for the next section, in which we will explain how these bounds can be propagated to a segment in the moduli space of triangles.
\subsection{First pass: separating the first 4 eigenvalues}
\label{section:first_pass}

We will focus on the lower bound of $\lambda_5$ needed to validate $\xi_{41}$, as the discussion for $\xi_{21}$ is analogous. In order to find a rigorous lower bound for the fifth eigenvalue of a triangle we will use a recent result found by Liu~\cite{Liu:framework-verified-eigenvalues}, which is similar to the one presented in~\cite{Carstensen-Gedicke:lower-bounds-eigenvalues} but simplifies the hypotheses and improves the constant involved. Both use the non-conforming Finite Element Method of Crouzeix--Raviart, which has the advantage over other conforming finite elements (such as \cite{Liu-Oishi:verified-eigenvalues-laplacian-polygons}) that the mass matrix is diagonal.

The Crouzeix--Raviart finite-element method uses a triangulation of the domain $\Omega$, which in our case we will take to be the trivial triangulation given by $N^2$ triangles with sides equal to $1/N$ of the original one. The basis functions are indexed by interior edges of the triangulation: if $E$ is a common edge of triangles $\tau$, $\tau'$, the basis function $\psi_E$ is the unique function supported on $\tau \cup \tau'$ such that restricted to each triangle is affine, takes the value $1$ in the midpoint of $E$ and the value $0$ in the midpoints of the other edges of $\tau$ and $\tau'$.

The weak formulation of the problem \eqref{laplacian-pde} reads:
\begin{align*}
\int_{\Omega} \nabla u \cdot \nabla \psi = \lambda \int_{\Omega} u \psi .
\end{align*}

Defining the coefficients of the stiffness and mass matrices $A = (a_{EF}), B = (b_{EF})$ by the bilinear forms
\[
  a_{EF} = \int_{\Omega} \nabla\psi_E\cdot\nabla\psi_F \qquad b_{EF} = \int_{\Omega} \psi_E\psi_F
\]

we are left to solve the discrete system
\begin{align}
\label{discrete-laplacian-pde-1}
A x = \lambda B x,
\end{align}

where the vector $x = (x_i)$ corresponds to the discrete solution $u_\text{d} = \sum x_i \psi_i$. For our choice of triangulation, $B$ is simply a multiple of the identity $2 |\Omega| I / \left(3N^2\right)$ (see~\cite[Remark 3.9]{Carstensen-Gedicke:lower-bounds-eigenvalues}), whereas $A$ is a sparse matrix whose entries involve in a straightforward way the geometry of the triangle. This allows us to invert $B$ in the system \eqref{discrete-laplacian-pde-1}, so we have to solve the following matrix eigenvalue problem
\begin{align}
\label{discrete-laplacian-pde-2}
M x = \lambda x, \quad M = 3N^2 A / (2 |\Omega|),
\end{align}

instead of the previous generalized one.
The main result that we will use is the following~\cite[Theorem 2.1 and Remark 2.2]{Liu:framework-verified-eigenvalues} linking solutions of the finite element problem \eqref{discrete-laplacian-pde-2} with solutions of the continuous one \eqref{laplacian-pde}:
\begin{thm}
  \label{liu}
  Consider a polygonal domain $\Omega$ with a triangulation so that each triangle has diameter at most $h$. Let $\lambda_{k}$ be the $k$-th solution of the eigenvalue problem \eqref{laplacian-pde} in $\Omega$ and $\lambda_{k,h}$ the $k$-th eigenvalue of the Crouzeix--Raviart discretized problem \eqref{discrete-laplacian-pde-2} in $\Omega$. Then
  \begin{equation}
    \label{eqn_liu}
    \frac{\lambda_{h,k}}{1+C_h^2 \lambda_{h,k}} \leq \lambda_k,
  \end{equation}
  where $C_h \leq 0.1893 h$ is a constant.
\end{thm}

In order to be able to deal with approximate eigenvalues we will need in addition the following a posteriori bound found in a lemma from~\cite[Theorem 15.9.1]{Parlett:symmetric-eigenvalue-book}, which given an approximate solution quantifies an upper bound on how far there has to be a true solution nearby.
\begin{lemma}
  \label{parlett}
  Let $(\tilde{\lambda}_h, \tilde{u}_h)$ be an approximate eigenpair solving \eqref{discrete-laplacian-pde-2} approximately such that $\tilde{\lambda}_h$ is closer to some $\lambda_h$ (which solves exactly \eqref{discrete-laplacian-pde-2}) than to any other solution of \eqref{discrete-laplacian-pde-2}. Suppose that the coefficient vector $\tilde{u}_h$ is normalised, $\|\tilde{u}_h\| = 1$. Then the (algebraic) residual $r := M\tilde{u}_h - \tilde{\lambda}_h \tilde{u}_h$ satisfies
  \[
    |\lambda_h - \tilde{\lambda}_h| \leq \|r\|.
  \]
\end{lemma}

\begin{rmk}
  \label{liu_error}
  We can combine Theorem~\ref{liu} with Lemma~\ref{parlett} thanks to the monotonicity of~\eqref{eqn_liu}, using $\lambda_{h,k} - \|r\|$ as a lower bound of $\lambda_{h,k}$ instead.
\end{rmk}

It is easy to obtain estimations $\tilde{\lambda}_h$ with a very small residual using standard eigenvalue routines. The hardest part here before applying the theorem is to check that they have indeed the correct index, i.e., that they are closer to the appropriate $\lambda_h$ than to any other discrete eigenvalue. This is why we need to control the whole spectrum of the discrete problem. We will do that by finding an approximately orthonormal basis which approximately diagonalizes the matrix, and proving that there is an exactly orthonormal basis nearby (with explicit bounds) of actual eigenvectors. This will give us useful bounds after application of Gershgorin disks theorem \cite{Gerschgorin:eigenvalues-theorem}. More precisely we will use the following:

\begin{lemma}
  \label{gram_schmidt}
  Let $v_1, \ldots, v_m$ be vectors in $\RR^m$ and $s > 0$ such that $|\langle v_i, v_j \rangle - \delta_{ij}| \leq s$, and suppose that $8ms < 1$. Then there exists an orthonormal set of vectors $w_1, \ldots, w_m \in \RR^m$ such that $\|v_i - w_i\| \leq \sqrt{3s}$\footnote{The constants in this lemma can be optimized, but for the sake of simplicity of the implementation, and since the errors in this part are very small, we have preferred to use simpler bounds.}.
  \begin{proof}
    Consider the Gram--Schmidt method applied to the vectors $v_1, \ldots, v_m$: we define recursively
    \[
      u_k := v_k - \sum_{i=1}^{k-1} u_i \frac{\langle u_i, v_k \rangle}{\|u_i\|^2}
    \]
    and the normalized vectors $w_k := u_k / \|u_k\|$. Let $\beta = s + 2 m s^2$. We will prove by induction that $\|u_k\|^2 \geq 1 - \beta$ and that $|\langle v_k, u_j \rangle| \leq \beta$ for all $j < k$. Indeed, for the first inequality,
    \begin{equation}
      \|u_k\|^2 = \langle v_k, v_k \rangle - \sum_{i=1}^{k-1} \frac{\langle v_k, u_i \rangle^2}{\|u_i\|^2} \geq 1 - s - m\frac{\beta^2}{1 - \beta}
      \label{close_iteration_norm}
    \end{equation}
    and for the second inequality,
    \begin{equation}
      |\langle v_k, u_j \rangle| \leq |\langle v_k, v_j\rangle| + \sum_{i=1}^{j-1} \frac{|\langle v_k, u_i\rangle \langle v_j, u_i\rangle|}{\|u_i\|^2} \leq s + m \frac{\beta^2}{1-\beta}.
      \label{close_iteration_product}
    \end{equation}
    We will be able to close the iteration as long as
    \[
      s + \frac{m \beta^2}{1 - \beta} \leq \beta
    \]
    (note that both~\eqref{close_iteration_norm} and~\eqref{close_iteration_product} give place to the same condition). Using the definition of $\beta$ this reduces to checking that
    $s(4m + 4m^2s + 2 + 4ms) \leq 1$, which holds true by virtue of the condition on $s$. Then we can bound the difference
    \[
      \|v_i - w_i\|^2 = \|v_i\|^2 + \|w_i\|^2 - 2 \langle v_i, w_i \rangle = \|v_i\|^2 + 1 - 2 \frac{\langle v_i, u_i \rangle}{\|u_i\|} = \|v_i\|^2 + 1 - 2 \|u_i\| \leq 2 + s - 2 \sqrt{1 - \beta}
    \]
which is bounded by $3s$.
  \end{proof}
\end{lemma}

We will use Lemma~\ref{gram_schmidt} in our proof in the following way: we will first find an approximate orthonormal basis of eigenvectors $\{v_i\}$, let $\tilde{Q}$ be the matrix that has them as columns, and compute an enclosure for the almost-diagonal matrix $\tilde{D} = \tilde{Q}^T M \tilde{Q}$. Let $Q$ be the matrix whose column vectors are given by the lemma, and construct $D = Q^T M Q$, which now has the same eigenvalues as $M$ due to the orthogonality of $Q$. We can obtain rigorous enclosures of the entries of $D$ using $\tilde{D}$ in this way:
\[
\begin{split}
       |D_{ij} - \tilde{D}_{ij}| &= |\langle w_i, M w_j \rangle - \langle v_i, M v_j \rangle| \\
  &\leq |\langle w_j - v_j, M v_i \rangle + \langle w_j - v_j, M (w_i - v_i) \rangle + \langle w_i - v_i, M v_j \rangle| \\
  &\leq \sqrt{3s} \left(\|M v_i \| + \|M v_j\|\right) + 4s \|M\|_2,
\end{split}
\]
where $s$ is as in the lemma and we have used the symmetry of $M$. Observe that the bounds for $\|M v_i\|$ can be computed explicitly, and that the upper bound $\|M\|_2 \leq \|M\|_\text{Frob}$ is easily computed.

Finally, applying Gershgorin's generalized disks theorem to the matrix $D$, of which we have sharp bounds, we can separate the spectrum of $M$ in at least two components, one of which will contain the first $5$ eigenvalues, and the other of which will contain the rest. This will allow us to rigorously identify the index of the $\lambda_{h,5}$, for which Lemma~\ref{parlett} provides a very sharp bound, by identifying $5$ disjoint enclosures of eigenvalues in the first component. Once we have this, application of Theorem~\ref{liu} and of Remark~\ref{liu_error} will provide a lower bound for $\lambda_5$.

\subsection{Second pass: getting narrow enclosures}

Our approach to find tight bounds for the eigenvalues of triangles uses the Method of Particular Solutions (MPS), introduced by Fox, Henrici and Moler in \cite{Fox-Henrici-Moler:approximations-bounds-eigenvalues} and more recently revived by Betcke and Trefethen \cite{Betcke-Trefethen:method-particular-solutions}. In this method, a function $u$ is written as a linear combination of functions $\phi_i$ ($1 \leq i \leq N$) that satisfy the equation $(\Delta + \lambda) \phi_i = 0$ on the full plane for a fixed $\lambda$. It is very important to remark that there are no boundary conditions on $\phi_i$ and the linear combination $u$ is chosen in such a way to ensure the boundary conditions of \eqref{laplacian-pde}. In other words, the coefficients are chosen to optimize the proximity of the function to the eigenspace of the actual eigenvalue $\lambda_j$. In order to achieve this, one tries to find a $\lambda$ and a non-zero function for which the boundary conditions are close to zero. This amounts to find a nontrivial solution $(c,\lambda)$ of
$$ A(\lambda) c = 0, $$

with $c_i$ being the coefficient vector such that $u(x) = \sum_{i=1}^N c_i \phi_i(x)$ and $A$ a rectangular matrix encoding the boundary conditions evaluated at certain points of the boundary $\partial \Omega$. Thus, the eigenvalues $\lambda$ will be the values for which the smallest singular value of $A$ is zero. In practice, $\lambda$ are found by using a golden ratio search on the smallest singular value of $A$. This provides a candidate $\lambda \in \RR$ and coefficients $c_i$ for which  $u$ can be computed with arbitrary precision.

The basis of functions $\phi_i$ that we will use was introduced very recently by Gopal and Trefethen~\cite{Gopal-Trefethen:new-laplace-solver-pnas,Gopal-Trefethen:laplace-solver-detailed} and offers root-exponential convergence. These functions are of two classes: the first are centered around points (that we will call charges) that accumulate exponentially near the vertices of the triangle along the outside of the angle bisector. More precisely, the $j$-th of them is located at a distance $\ell e^{-\sigma j / \sqrt{n_c}}$ of the vertex, where the parameters we have used are $\ell = 1.0, \sigma = 2.5$. Here $n_c$ is the number of charges and $0 \leq j < n_c$; we have used $n_c = 7$. Using polar coordinates around the charges, the three functions centered around every charge take the form
\[
  \phi_\text{ext}(r, \theta) = Y_0\left(r \sqrt{\lambda}\right), \qquad
  \phi_\text{ext}^\text{c}(r, \theta) = Y_1\left(r \sqrt{\lambda}\right) \cos \theta, \qquad
  \phi_\text{ext}^\text{s}(r, \theta) = Y_1\left(r \sqrt{\lambda}\right) \sin \theta.
\]

The functions of the second class are all centered around the same point, which we have chosen to be the centroid of the triangle. Using polar coordinates around it, they take the form
\[
  \phi_0(r, \theta) = J_0\left(r \sqrt{\lambda}\right), \qquad
  \phi_j^\text{c}(r, \theta) = J_j\left(r \sqrt{\lambda}\right) \cos j\theta, \qquad
  \phi_j^\text{s}(r, \theta) = J_j\left(r \sqrt{\lambda}\right) \sin j\theta.
\]
for $1 \leq j \leq d$; we have used $d = 10$. Here $J_j$ and $Y_j$ are the Bessel functions of the first and second kind, respectively. Note that the first class of functions are singular at the charge point, which allows us to develop the correct behavior at the vertex of the triangle, where the eigenfunction is singular. The second class of functions allows us to fit the boundary values near the vertices.

This method produces very accurate eigenvalues and eigenfunctions but with no guaranteed bounds. In order to obtain these, we will make use of the method developed by Barnett and Hassell \cite{Barnett-Hassell:quasi-ortogonality-dirichlet-eigenvalues}, that bounds the error if we can control the $L^2$ norm of the candidate eigenfunction on the border. In addition to an improvement in the constant by a factor of $\sqrt{\lambda}$, which makes it optimal in the high eigenvalue limit, the $L^2$ control is a lot more efficient in practice than the $L^\infty$ bounds originally developed in~\cite{Fox-Henrici-Moler:approximations-bounds-eigenvalues,Moler-Payne:bounds-eigenvalues}, because the approximate eigenfunctions generated with typical MPS bases tend to be very irregular around the vertices but regular and small in most of the boundary. Their method is however optimized for high-index eigenvalues, so we adapt some of the steps to our case of small eigenvalues, simplifying the bookkeeping of the explicit constants.

We summarize the main results that we will use. Let $\Omega$ be a triangle and $u$ a nonzero smooth function on $\Omega$ such that $(\Delta + \lambda) u = 0$. We will measure how good $u$ is as an approximation to an eigenfunction using the tension
\begin{equation}
  \label{tension_def}
  t[u] := \frac{\|u\|_{L^2(\partial \Omega)}}{\|u\|_{L^2(\Omega)}}.
\end{equation}
Let $\lambda_j, u_j$ be the sequence of eigenvalues and eigenfunctions of $\Omega$, satisfying $(\Delta + \lambda_j) u_j = 0$ with zero Dirichlet boundary conditions. Let $v_j$ be the normal derivative of $u_j$, defined on the regular part of $\partial \Omega$. We define the operator
\[
  A(\lambda) = \sum_{\lambda_j} \frac{v_j \langle v_j, \cdot \rangle}{(\lambda - \lambda_j)^2},
\]
and decompose it as a sum of three:
\[
  A_\text{near}(\lambda) = \sum_{|\lambda - \lambda_j| \leq \sqrt{\lambda}} \frac{v_j \langle v_j, \cdot \rangle}{(\lambda - \lambda_j)^2},
\]
\[
  A_\text{far}(\lambda) = \sum_{\lambda / 2 \leq \lambda_j \leq 2 \lambda, |\lambda - \lambda_j| > \sqrt{\lambda}} \frac{v_j \langle v_j, \cdot \rangle}{(\lambda - \lambda_j)^2},
\]
\[
  A_\text{tail}(\lambda) = \sum_{\lambda_j < \lambda / 2 \text{ or } \lambda_j > 2 \lambda} \frac{v_j \langle v_j, \cdot \rangle}{(\lambda - \lambda_j)^2}
\]
where $\langle \cdot, \cdot \rangle$ is the standard inner product on $L^2(\partial \Omega)$.

This operator is useful because its norm is controlled from below by the tension (see~\cite[Section 3]{Barnett-Hassell:quasi-ortogonality-dirichlet-eigenvalues}):
\begin{equation}
  \label{tension}
  t[u]^{-2} \leq \|A(\lambda)\|.
\end{equation}

Moreover we have the following bounds from \cite[Lemma 4.1]{Barnett-Hassell:quasi-ortogonality-dirichlet-eigenvalues} and \cite[Lemma 4.2]{Barnett-Hassell:quasi-ortogonality-dirichlet-eigenvalues}:
\begin{equation}
  \label{bound_far}
  \|A_\text{far}(\lambda)\| \leq C_1,
\end{equation}
\begin{equation}
  \label{bound_tail}
  \|A_\text{tail}(\lambda)\| \leq C_2 \lambda^{-1/2},
\end{equation}
that hold for all $\lambda > 1$\footnote{This condition obviously holds for our triangles, for example by the Faber--Krahn inequality.}, with explicit constants $C_1, C_2$ given below. For the near term, since we are working with very low eigenvalues, $\sqrt{\lambda}$ is actually small enough that only the summand with $\lambda_j = \lambda$ appears:
\begin{equation}
  \label{near_term}
  \|A_\text{near}(\lambda)\| = \frac{\|v_j\|_{L^2(\partial \Omega)}^2}{(\lambda - \lambda_j)^2}.
\end{equation}

For convex domains, like in our case, Section~6 of~\cite{Barnett-Hassell:quasi-ortogonality-dirichlet-eigenvalues} offers explicit bounds for the constants. In particular, they are based on another constant $C_\Omega$ coming from a boundary quasi-orthogonality inequality (see~\cite[Lemma 2.1]{Barnett-Hassell:quasi-ortogonality-dirichlet-eigenvalues}) that we will not discuss in more detail. We will simply use the bounds $C_1 \leq 2 \pi^2 C_\Omega / 3 < 7 C_\Omega$ and $C_2 < 7 C_\Omega$ quoted from the end of page 1058 of the paper.

$C_\Omega$ can be taken as any constant that makes Lemma 2.1 in \cite{Barnett-Hassell:quasi-ortogonality-dirichlet-eigenvalues} hold; in particular, in view of the proof of the lemma and since $\lambda > 1$, one can take $C_\Omega = 4(C_a + C_a') + \sqrt{2} C_a''$, where these three newly introduced constants depend on a vector field on $\Omega$ with certain properties. Luckily, for strictly star-shaped domains, one can make a trivial choice of such a vector field and obtain simple expressions for the constants. The concrete expressions can be quoted from the first bullet point of~\cite[Section 6]{Barnett-Hassell:quasi-ortogonality-dirichlet-eigenvalues}:
\[
  C_a = \frac{\sup_{\partial \Omega} (\mathbf{x} \cdot \mathbf{n})}{\inf_{\partial \Omega} (\mathbf{x} \cdot \mathbf{n})}, \qquad
  C_a' = \frac{1}{\inf_{\partial \Omega} (\mathbf{x} \cdot \mathbf{n})}, \qquad
  C_a'' = 0,
\]
where $\mathbf{n}$ is the outer normal vector field in the regular part of the boundary. Since the problem is coordinate-invariant, we are free to chose an origin, and the optimal choice is the incenter of the triangle. In this case, both the supremum and the infimum are equal to the inradius $\rho$, so we can take $C_\Omega = 4 (1+\rho) / \rho$.

Putting \eqref{tension}-\eqref{near_term} together we have
\begin{equation}
  \label{tension_bound}
  t[u]^{-2} \leq \frac{\|v_j\|_{L^2(\partial \Omega)}^2}{(\lambda - \lambda_j)^2} + 7 C_\Omega (1 + \lambda^{-1/2}).
\end{equation}

Finally, recall Rellich's formula~\cite{Rellich:darstellung-eigenwerte}:
\[
  \int_{\partial \Omega} (\partial_n u_j)^2 (\mathbf{x} \cdot \mathbf{n}) \mathrm{d}s = 2 \lambda_j.
\]
For our choice of origin of coordinates, this just gives us $\|v_j\|_{L^2(\partial \Omega)}^2 = 2 \lambda_j / \rho$. Inserting this into~\eqref{tension_bound} and rearranging we have proved:

\begin{prop}
  \label{second_pass}
  The distance $d$ from $\lambda$ to the spectrum of the Laplacian on $\Omega$, can be bounded above by
  \[
    d \leq \sqrt{\frac{2 \tilde{\lambda_j}}{\rho t[u]^{-2} - 28 (1 + \rho) (1 + \lambda^{-1/2})}}
  \]
  where $\tilde{\lambda_j}$ is an upper bound for $\lambda_j$\footnote{Obtaining this upper bound is a minor problem: one can use first a trivial one, like the one obtained in the first pass, to obtain a first enclosure and then use it to bootstrap and get a better one.}.
\end{prop}

\section{Proof of Theorem \ref{main_thm} and reduction to a finite number of triangles}
\label{section:twotriangles}

In this section we will give specific details of the triangles $T_\mathrm{A}$ and $T_\mathrm{B}$ that we construct in Theorem \ref{main_thm}. We first notice that because of the scaling of the problem, one can fix two vertices at $(0,0)$ and $(1,0)$ and require the quotients $\xi_{21} = \lambda_2 / \lambda_1$ and $\xi_{41} = \lambda_4 / \lambda_1$ to be equal. We will parametrize the search space by means of the coordinates $(c_x,c_y)$ of the third vertex of the triangle. 

Let $A = (0.63500, 0.27500)$ and $B = (0.84906, 0.31995)$ and consider the parallelograms in this space with the following sets of vertices:
\begin{itemize}
  \item $A + v_{21,A} + v_{41,A}$, $A + v_{21,A} - v_{41,A}$, $A - v_{21,A} - v_{41,A}$, $A - v_{21,A} + v_{41,A}$;
  \item $B + v_{21,B} + v_{41,B}$, $B + v_{21,B} - v_{41,B}$, $B - v_{21,B} - v_{41,B}$, $B - v_{21,B} + v_{41,B}$;
\end{itemize}

where the explicit values of the vectors are
\begin{itemize}
  \item $v_{21,A} = ( 0.004610608896618232,   0.0012403688839389946)$
  \item $v_{41,A} = (-0.0041659682109460045, -0.000511581170421992)$
  \item $v_{21,B} = ( 0.0028159587453638808,  0.0020776257941965285)$
  \item $v_{41,B} = ( 0.007180726583099708,   0.00213029677112299)$
\end{itemize}

\begin{figure}
  \centering
  \input{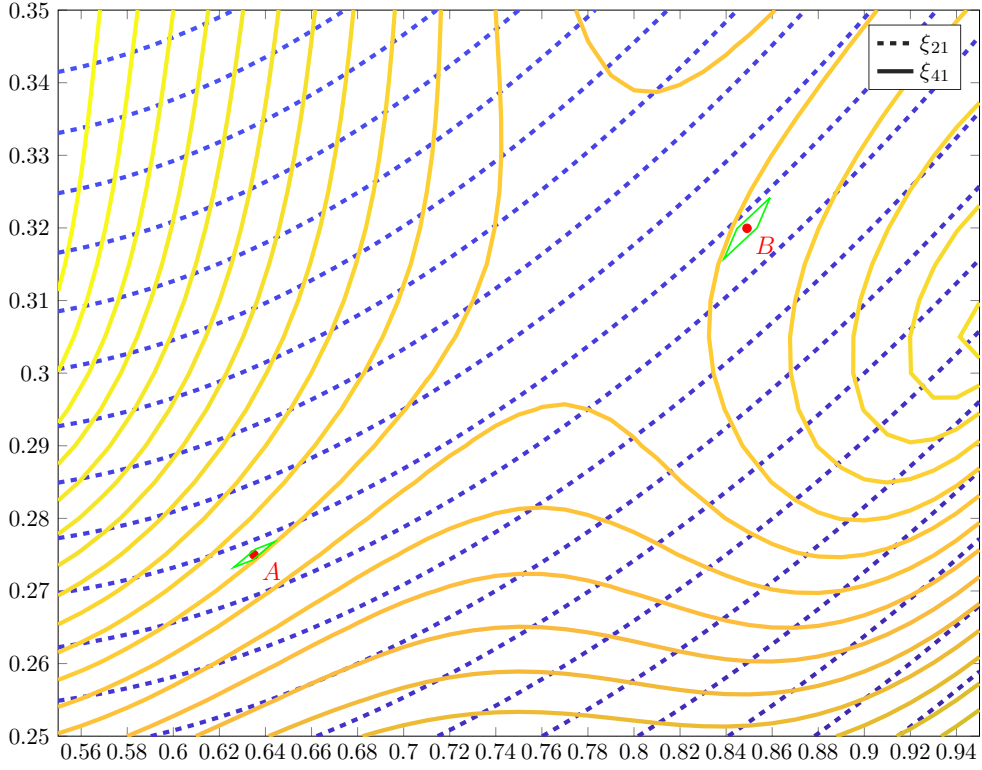}
  \caption{(Numerical) plot of the level sets of the quotients $\xi_{21}$ (discontinuous lines) and $\xi_{41}$ (continuous lines) around the region of interest. The points $A$ and $B$ are shown in red and the validated parallelograms in green.}
  \label{tworegions}
\end{figure}

These are drawn in Figure~\ref{tworegions} together with a numerical contour plot of the eigenvalue quotients. The vectors have been chosen with the heuristic that moving along them keeps one of $\xi_{21}, \xi_{41}$ approximately constant, while incrementing the other by a fixed amount. In particular, their lengths have been adjusted to optimize the total time of the calculations, as shortening one pair of sides of the parallelogram reduces the total time to validate them, but it also reduces the margin $\xi - \bar{\xi}$ in the validation of the other pair of sides, thereby increasing their validation time.

We can reduce Theorem \ref{main_thm} to a computation on the boundary of the parallelograms by the following proposition:

\begin{prop}
 Let the parallelograms be defined as before and consider $\xi_{21} - \bar{\xi}_{21}$ and $\xi_{41} - \bar{\xi}_{41}$, with $\bar{\xi}_{21}:=1.67675$ and $\bar{\xi}_{41}:=2.99372$. If the following condition is satisfied:
 \[
   \parbox{0.8\textwidth}{each of $\xi_{21} - \bar{\xi}_{21}$ and $\xi_{41} - \bar{\xi}_{41}$ has a constant and opposite sign in opposite edges of the parallelograms} \tag{C}
 \]
then there exist two triangles $T_\mathrm{A}$ and $T_\mathrm{B}$ such that 

$$\xi_{21}(T_\mathrm{A}) = \xi_{21}(T_\mathrm{B}) = \bar{\xi}_{21} \text{ and } \xi_{41}(T_\mathrm{A}) = \xi_{41}(T_\mathrm{B}) = \bar{\xi}_{41}.$$
\end{prop}
 The proof follows from the continuity of eigenvalues and the Poincar\'e--Miranda theorem \cite{Miranda:PoincareMiranda}, which we recall here for completeness:

\begin{thm}
  Given two continuous functions $f, g : [-1, 1]^2 \to \RR$ such that $f(x, y)$ has the same sign as $x$ when $x = \pm 1$ and $g(x, y)$ has the same sign as $y$ when $y = \pm 1$, there exists a point $(x, y) \in [-1, 1]^2$ such that $f(x, y) = g(x, y) = 0$.
\end{thm}

Note that the slightly more general result of Remark~\ref{rmk_open} follows from the fact that with the computer we will check only a finite number of strict inequalities, so they will remain true if the values of $\bar{\xi}_{21}$ and $\bar{\xi}_{41}$ are perturbed by a small amount.

Our goal now is to check the aforementioned sign condition on the sides of the parallelograms. Instead of working directly with the full side, we will show that it is enough to check condition (C) on a finite set of triangles through stability estimates. These will follow from the monotonicity of the eigenvalues with respect to the inclusion of domains and the scaling properties of the eigenvalues.

\begin{figure}[h!]
\centering
\includegraphics[scale=0.5]{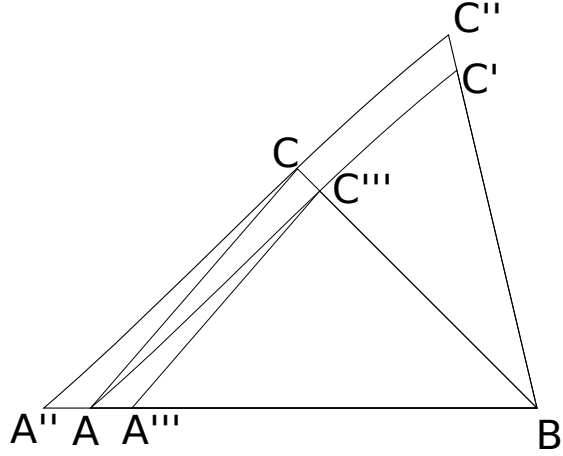}
    \caption{Two pairs of homothetic triangles:  $ABC$ and $A'''BC'''$, and $ABC'$ and $A''BC''$.}
  \label{fig_homothety}
\end{figure}

\begin{lemma}
  \label{homothety1}
  Let $T$ and $T'$ be two triangles, whose vertices are $A = (0, 0)$, $B = (1, 0)$, and $C = (c_x, c_y), C' = (c_x', c_y')$ respectively ($c_y, c_y' > 0$). Consider the cross products $p = \overrightarrow{AC} \times \overrightarrow{AC'} = c_x c_y' - c_y c_x'$ and $q = \overrightarrow{BC} \times \overrightarrow{BC'} = (c_x - 1) c_y' - c_y (c_x' - 1)$. Then:
  \begin{itemize}
    \item If both $p, q < 0$, there is a homothety of $T'$ by a factor $1 - p / c_y'$ that contains $T$.
    \item If both $p, q > 0$, there is a homothety of $T'$ by a factor $1 + q / c_y'$ that contains $T$.
  \end{itemize}
  \begin{proof}
    The proof is very similar in the two cases, so we will only do it for the first one. We want to find the homothety of scale $1 + r$ that keeps the vertex $B$ of triangle $T'$ fixed and such that the image of its opposite side contains vertex $C$ of $T$ (see Figure \ref{fig_homothety} for a picture of the homothetic triangles $ABC$ and $A'''BC'''$, and $ABC'$ and $A''BC''$ respectively). The condition becomes simpler once we apply an inverse homothety to $T$ and $T'$, so that it results in the points
    \[
      A, C''' = \frac{C + r B}{1 + r}, C'
    \]
    being aligned. The solution is $r = -p / c_y'$, which is positive by our condition. Moreover, triangle $T$ lies below this homothety of $T'$ because the vectors $\overrightarrow{BC}$ and $\overrightarrow{BC'}$ are in the correct orientation due to the condition $q < 0$. This suffices to check that $T$ is contained in this homothety.
  \end{proof}
\end{lemma}

\begin{lemma}
  \label{homothety2}
  With the same notation as in Lemma~\ref{homothety1},
  \begin{itemize}
    \item If $p > 0$ and $q < 0$, then $T \subset T'$.
    \item If $p < 0$ and $q > 0$, there is a homothety of $T'$ by a factor $c_y / c_y' > 1$ that contains $T$.
  \end{itemize}
  \begin{proof}
    In the first case, the conditions on the signs of the cross products of the side vectors is equivalent to $T$ being contained in $T'$. In the second case, the relative orientations of the sides guarantee that a homothetic triangle to $T'$ of the same height as $T$ whose top vertex coincides with $C$ will contain $T$, and the ratio of this homothety is clearly $c_y / c_y'$.
  \end{proof}
\end{lemma}

Using two reversed inclusions from the previous lemmas we can prove:

\begin{lemma}
  \label{stability}
  Let $T$ be a triangle as above, and consider perturbations of the third vertex of the form $C + t v$ defining triangles $T^{(t)}$, for $t \in [-\ell, \ell]$, where $v = (v_x, v_y)$. Let $\lambda_n, \lambda_n^{(t)}$ be the $n$-th Dirichlet eigenvalues of triangles $T, T^{(t)}$, respectively, and define $\xi_{n1}^{(t)}$ as the obvious eigenvalue quotient. Then we distinguish two cases depending on $p_v = \overrightarrow{AC} \times v$ and $q_v = \overrightarrow{BC} \times v$:
  \begin{itemize}
    \item If $p_v$ and $q_v$ both have the same sign, then for all $t \in [-\ell, \ell]$
  \[
    |\xi_{n1}^{(t)} - \xi_{n1}| \leq \xi_{n1} \left[\left(1 + \ell \frac{|p_v|}{c_y - \ell|v_y|}\right)^2 \left(1 + \ell \frac{|q_v|}{c_y - \ell|v_y|}\right)^2 - 1\right].
  \]
    \item If $p_v$ and $q_v$ have different signs, then for all $t \in [-\ell, \ell]$
  \[
    |\xi_{n1}^{(t)} - \xi_{n1}| \leq \xi_{n1} \left[\left(\frac{c_y}{c_y - \ell |v_y|}\right)^2 - 1\right].
  \]
  \end{itemize}

  \begin{proof}
    For the first part, note that either of the conditions of Lemma~\ref{homothety1} holds for $T$ and $T^{(t)}$ whenever $t \neq 0$: indeed, $p = \overrightarrow{AC} \times (\overrightarrow{AC} + tv) = t p_v$ and $q = \overrightarrow{BC} \times (\overrightarrow{BC} + tv) = t q_v$ have the same sign if $p_v$ and $q_v$ do. Note also that reversing the roles of $T$ and $T^{(t)}$ causes both $p$ and $q$ to change sign and hence the Lemma still applies, but with the other condition. Combining this with eigenvalue monotonicity ($\Omega \subseteq \Omega' \Rightarrow \lambda_n(\Omega) \geq \lambda_n(\Omega')$) and scaling ($\lambda_n(\rho \Omega) = \rho^{-2} \lambda_n(\Omega)$) yields either
    \[
      \begin{split}
      \lambda_n \geq \left(1 + \frac{|t p_v|}{c_y + t v_y}\right)^{-2} \lambda_n^{(t)} \geq \left(1 + \ell \frac{|p_v|}{c_y - \ell |v_y|}\right)^{-2} \lambda_n^{(t)}, \\
      \lambda_n^{(t)} \geq \left(1 + \frac{|t q_v|}{c_y}\right)^{-2} \lambda_n > \left(1 + \ell \frac{|q_v|}{c_y - \ell |v_y|}\right)^{-2} \lambda_n 
      \end{split}
    \]
    or
    \[
      \begin{split}
      \lambda_n \geq \left(1 + \frac{|t q_v|}{c_y + t v_y}\right)^{-2} \lambda_n^{(t)} \geq \left(1 + \ell \frac{|q_v|}{c_y - \ell |v_y|}\right)^{-2} \lambda_n^{(t)}, \\
      \lambda_n^{(t)} \geq \left(1 + \frac{|t p_v|}{c_y}\right)^{-2} \lambda_n > \left(1 + \ell \frac{|p_v|}{c_y - \ell |v_y|}\right)^{-2} \lambda_n
      \end{split}
    \]
    independently of $n$. In either case the conclusion is that
    \[
      \left(1 + \ell \frac{|p_v|}{c_y - \ell |v_y|}\right)^{-2} \left(1 + \ell \frac{|q_v|}{c_y - \ell |v_y|}\right)^{-2} \leq \frac{\lambda_n^{(t)} / \lambda_n}{\lambda_1^{(t)} / \lambda_1} = \frac{\xi_{n1}^{(t)}}{\xi_{n1}} \leq \left(1 + \ell \frac{|p_v|}{c_y - \ell |v_y|}\right)^2 \left(1 + \ell \frac{|q_v|}{c_y - \ell |v_y|}\right)^2
    \]
    and from here the first claim readily follows. For the second one the argument is similar: either
    \[
      \left(\frac{c_y - \ell |v_y|}{c_y}\right)^{-2} \lambda_n^{(t)} \geq \left(\frac{c_y + \ell |v_y|}{c_y}\right)^2 \lambda_n^{(t)} \geq \left(\frac{c_y'}{c_y}\right)^2 \lambda_n^{(t)} \geq \lambda_n \geq \lambda_n^{(t)},
    \]
    or 
    \[
      \left(\frac{c_y - \ell |v_y|}{c_y}\right)^2 \lambda_n^{(t)} \leq \left(\frac{c_y'}{c_y}\right)^2 \lambda_n^{(t)} \leq \lambda_n \leq \lambda_n^{(t)},
    \]
    independently of $n$ too. Similarly, in either case the conclusion is that
    \[
      \left(\frac{c_y - \ell |v_y|}{c_y}\right)^2 \leq \frac{\lambda_n^{(t)} / \lambda_n}{\lambda_1^{(t)} / \lambda_1} = \frac{\xi_{n1}^{(t)}}{\xi_{n1}} \leq \left(\frac{c_y - \ell |v_y|}{c_y}\right)^{-2}
    \]
    and the claim follows.
  \end{proof}
\end{lemma}

In the last section we discuss how to pick the triangles for which we check Condition (C), and how the two passes of Section~\ref{section:onetriangle} are implemented, thus concluding the proof.

\section{Details of the implementation}
\label{section:implementation}
The computer-assisted part of the proof is done with three main programs that take care of independent parts of the validation. All the rigorous computations are performed using the validated arithmetics library Arb, developed by Fredrik Johansson~\cite{Johansson:Arb}, which can be found at \url{http://arblib.org}. The three programs are:
\begin{enumerate}[label=(\alph*)]
  \item \label{item:valxi} The program implemented in \verb|valxi_main.cc| takes as input data specifying a side of one of the two quadrilaterals, a total number $N$ of subdivisions of the side, and an integer $1 \leq c \leq N$ indexing the subdivision to validate. Then it computes narrow enclosures for two eigenvalues (which program~\ref{item:position} will check a posteriori that correspond to $\lambda_1$ and $\lambda_k$, for $k = 2$ or $4$), of the triangle parametrized by the center of the segment. Finally it outputs them together with an enclosure for $\xi_{k1}$ which is valid for the whole subdivided segment.
  \item \label{item:interm} \verb|interm_main.cc| reads data specifying one of the sides in which we are validating $\xi_{41}$, and uses the same machinery as~\ref{item:valxi} to find enclosures for two eigenvalues (which program~\ref{item:position} will later confirm to be $\lambda_2$ and $\lambda_3$) of the triangle corresponding to the center of the side. Then it outputs bounds propagated to the whole side (it turns out that here with $N=1$ segment we get sufficiently small errors).
  \item \label{item:position} Finally, \verb|position_main.cc| is responsible for the first pass described above: separating the first $k = 2$ or $4$ eigenvalues from the rest of the spectrum and checking that the eigenvalues obtained in~\ref{item:valxi} do have the claimed indices. This is done by first finding a rigorous lower bound for $\lambda_{k+1}$ at the center of the side using the FEM technique, and propagating it to the whole side. In the case $k=2$, it is enough to check that the eigenvalues found by~\ref{item:valxi} on each subdivided interval lie in disjoint enclosures in $(0, \lambda_3)$. In the case $k=4$, we need to check that the two eigenvalues given by~\ref{item:valxi} and the two intermediate eigenvalues given by~\ref{item:interm} give pairwise disjoint enclosures in $(0, \lambda_5)$ for each subdivided interval.
\end{enumerate}

In our implementation, for part~\ref{item:valxi} we have divided each side into $N = 40$ smaller segments. For~\ref{item:interm} and~\ref{item:position} we have not needed any subdivision. Therefore the bulk of the time is taken by the first program: one run takes approximately $8$ minutes, so the total running time is of approximately $42$ hours. This can be of course run in parallel, in our case reducing the time to less than one hour with $80$ machines.

In programs~\ref{item:valxi} and \ref{item:interm}, eigenvalues are searched around the approximate numeric eigenvalues of a FEM discretization of the problem, and then they are refined and rigorously validated using the techniques explained in the second pass section. These approximate eigenvalues are found using the numeric linear algebra library ALGLIB. The refined search for the eigenvalues is done by minimizing the least singular value of a matrix involving boundary values of the basis functions, as described in \cite{Betcke-Trefethen:method-particular-solutions}; we have chosen $300$ boundary points per side, distributed as Chebyshev nodes, and $40$ interior points randomly chosen with a fixed seed. The singular value is computed using the ALGLIB library again, and the minimum is found by a golden ratio search.

\begin{figure}[h!]
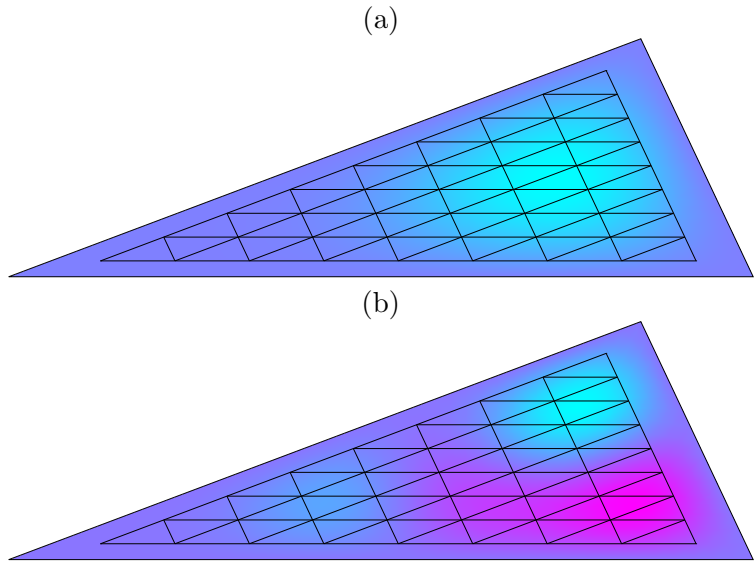

  \centering
    \vbox{
      (a)

      \resizebox{0.6\textwidth}{!}{\input{first_eigenfunction_B.tex}}

      (b)
      
      \resizebox{0.6\textwidth}{!}{\input{fourth_eigenfunction_B.tex}}
    }
    \caption{Grid used to validate a lower bound for $\|u\|_{L^2(\Omega)}$ for triangle $T_\mathrm{B}$, with (a) the first eigenfunction and (b) the second eigenfunction plotted on top.}
  \label{grid_plot}
\end{figure}

For parts~\ref{item:valxi} and~\ref{item:interm}, in order to apply Proposition~\ref{second_pass} we need an upper bound on the tension defined in~\eqref{tension_def}, which amounts to finding an upper bound for the $L^2$ norm of $u$ on $\partial \Omega$ and a lower bound on its interior $L^2$ norm. We now discuss the details of both:

\begin{itemize}

\item Upper bound of $\|u\|_{L^2(\partial \Omega)}$: 

For the upper bound on the boundary norm, the sides of the triangle are divided into small segments, starting with a Chebyshev node distribution and subdividing into two smaller segments whenever the supremum of $u^2$ in the interval is larger than a fixed threshold. For each of these segments $\sigma$, linearly parametrized by $p : [-1, 1] \to \sigma \subset \partial \Omega$, the integral of $v(t) := u^2(p(t))$ is evaluated by using a quadrature for the corresponding Taylor expansion, and adding a computed upper bound of the error term. Note that to control both the coefficients of the expansion and the error term, the derivatives of $v$ need to be evaluated both in the midpoint of the segment and in the whole segment using interval arithmetic.

More precisely, if $v(t) = \sum_{j = 0}^m v_j t^j + r(t)$, where $v_j = \frac{v^{(j)}(0)}{j!}$ are the Taylor coefficients and the residue term $r(t)$ can be bounded by $R = \sup_{t \in [-1, 1]} \left| \frac{v^{(m+1)}(t)}{(m+1)!} \right|$, then we can use the bound
\[
  \int_\sigma u^2 \leq |\sigma| \left(\sum_{i = 0}^{\lfloor m/2 \rfloor} \frac{v_{2i}}{2i + 1} + \frac{R}{m + 2} \right)
\]
assuming $m$ is odd.

\item Lower bound of $\|u\|_{L^2(\Omega)}$: 

For the interior norm lower bound, we will only bound from below the integral of $u^2$ in an interior region of $\Omega$, because the contributions to the $L^2$ norm near the boundary are small due to the Dirichlet condition. We will do this by partitioning this interior region into a triangular grid and using a form of the minimum principle to deduce a lower bound of $u^2$ in each triangle $\tau$ of the grid from a lower bound of $u^2$ on $\partial \tau$. More precisely, suppose $|u| \geq b$ on $\partial \tau$, and assume without loss of generality that $u > 0$ on $\partial \tau$. If $u < 0$ anywhere inside $\tau$, the region $U$ on which $u$ is negative has $\lambda$ as its first Dirichlet eigenvalue, so by the Faber--Krahn inequality, $|U|$ is bounded below. However, $U \subset \tau$ implies $|U| \leq |\tau|$, so for a small enough partition this will be a contradiction. Therefore $-\Delta u = \lambda u > 0$ on $\tau$, and by the minimum principle, $u \geq \sup_{\partial \tau} u \geq b$ on $\tau$.

If $u$ changes sign on the boundary, the lower bound of $|u|$ that we will get will be $0$, so we will be ignoring that triangle. It turns out that this strategy is enough to get a fraction of $\|u\|^2$ as a lower bound.

In our case we have used a grid of $8^2$ triangles that partitions an interior triangle of side length $0.8$ times the original one (see Figure~\ref{grid_plot}). On the perimeter of the grid, lower bounds of $u^2$ are obtained by subdividing the segments into halves until $u^2$ is computed on each segment with a precision smaller than a given threshold. Evaluation of $u^2$ on a segment $\sigma$ is done by means of the Taylor expansion again, controlling both the derivatives of $v(t) = u^2(p(t))$ evaluated at $t = 0$ and at the whole interval $t \in [-1, 1]$ by means of interval arithmetic.
\end{itemize}

Due to the oscillatory behavior of the solutions, especially near the vertices, and the large amount of cancellation between the summands, a high order Taylor expansion (order $m = 25$) has been used to evaluate $u^2$. Since the derivatives of the Bessel functions, which take most of the time of the computation, can be computed recursively at each point, it is more efficient to use a high order expansion and fewer intervals to speed up the validation.
  
For part~\ref{item:position}, the approximate orthonormal basis $\{v_i\}$ described at the end of section~\ref{section:first_pass} is constructed using the eigendecomposition method for symmetric matrices from the numerical library ALGLIB. The rest of the calculations were easily performed using Arb.

\section*{Acknowledgements}

We are grateful to Princeton University, where this research was performed, and to the VSRC Program and the Department of Mathematics for partially supporting the second author's stay. JGS was partially supported by NSF through Grant NSF DMS-1763356. GO was partially supported by the CFIS Mobility Grant, by the MOBINT-MIF Scholarship from AGAUR, and by the Severo Ochoa grant program from ICMAT. We also thank Princeton University for computing facilities (Polar cluster).

\bibliographystyle{abbrv}
\bibliography{references}

\def\cprime{$'$}
\begin{thebibliography}{10}

\bibitem{Antunes-Freitas:inverse-spectral-problem-triangles}
P.~R.~S. Antunes and P.~Freitas.
\newblock On the inverse spectral problem for {E}uclidean triangles.
\newblock {\em Proc. R. Soc. Lond. Ser. A Math. Phys. Eng. Sci.},
  467(2130):1546--1562, 2011.

\bibitem{Antunes-Valtchev:mfs-corners-cracks}
P.~R.~S. Antunes and S.~S. Valtchev.
\newblock A meshfree numerical method for acoustic wave propagation problems in
  planar domains with corners and cracks.
\newblock {\em J. Comput. Appl. Math.}, 234(9):2646--2662, 2010.

\bibitem{Ashbaugh-Benguria:ppw-conjecture}
M.~S. Ashbaugh and R.~D. Benguria.
\newblock A sharp bound for the ratio of the first two eigenvalues of
  {D}irichlet {L}aplacians and extensions.
\newblock {\em Annals of Mathematics}, 135(3):601--628, 1992.

\bibitem{Barnett-Hassell:quasi-ortogonality-dirichlet-eigenvalues}
A.~H. Barnett and A.~Hassell.
\newblock Boundary quasi-orthogonality and sharp inclusion bounds for large
  {D}irichlet eigenvalues.
\newblock {\em SIAM J. Numer. Anal.}, 49(3):1046--1063, 2011.

\bibitem{Beattie-Goerisch:lower-bounds-eigenvalues}
C.~Beattie and F.~Goerisch.
\newblock Methods for computing lower bounds to eigenvalues of self-adjoint
  operators.
\newblock {\em Numer. Math.}, 72(2):143--172, 1995.

\bibitem{Behnke-Goerisch:inclusions-eigenvalues}
H.~Behnke and F.~Goerisch.
\newblock Inclusions for eigenvalues of selfadjoint problems.
\newblock In {\em Topics in validated computations ({O}ldenburg, 1993)},
  volume~5 of {\em Stud. Comput. Math.}, pages 277--322. North-Holland,
  Amsterdam, 1994.

\bibitem{Betcke:generalized-svd-mps}
T.~Betcke.
\newblock The generalized singular value decomposition and the method of
  particular solutions.
\newblock {\em SIAM J. Sci. Comput.}, 30(3):1278--1295, 2008.

\bibitem{Betcke-Trefethen:method-particular-solutions}
T.~Betcke and L.~N. Trefethen.
\newblock Reviving the method of particular solutions.
\newblock {\em SIAM Rev.}, 47(3):469--491, 2005.

\bibitem{Carstensen-Gedicke:lower-bounds-eigenvalues}
C.~Carstensen and J.~Gedicke.
\newblock Guaranteed lower bounds for eigenvalues.
\newblock {\em Math. Comp.}, 83(290):2605--2629, 2014.

\bibitem{Chang-DeTurck:n-eigenvalues-triangle}
P.-K. Chang and D.~DeTurck.
\newblock On hearing the shape of a triangle.
\newblock {\em Proc. Amer. Math. Soc.}, 105(4):1033--1038, 1989.

\bibitem{Croke-Sharafutdinov:spectral-rigidity-negatively-curved-manifold}
C.~B. Croke and V.~A. Sharafutdinov.
\newblock Spectral rigidity of a compact negatively curved manifold.
\newblock {\em Topology}, 37(6):1265--1273, 1998.

\bibitem{Durso:phd-thesis}
C.~Durso.
\newblock {\em On the inverse spectral problem for polygonal domains}.
\newblock ProQuest LLC, Ann Arbor, MI, 1988.
\newblock Thesis (Ph.D.)--Massachusetts Institute of Technology.

\bibitem{Enciso-GomezSerrano:spectral-semiregular}
A.~Enciso and J.~G\'omez-Serrano.
\newblock Spectral determination of semi-regular polygons.
\newblock {\em Arxiv preprint arXiv:1709.05960}, 2017.

\bibitem{Fairweather-Karageorghis:mfs-survey}
G.~Fairweather and A.~Karageorghis.
\newblock The method of fundamental solutions for elliptic boundary value
  problems.
\newblock volume~9, pages 69--95. 1998.
\newblock Numerical treatment of boundary integral equations.

\bibitem{Fox-Henrici-Moler:approximations-bounds-eigenvalues}
L.~Fox, P.~Henrici, and C.~Moler.
\newblock Approximations and bounds for eigenvalues of elliptic operators.
\newblock {\em SIAM J. Numer. Anal.}, 4:89--102, 1967.

\bibitem{Gerschgorin:eigenvalues-theorem}
S.~A. Gershgorin.
\newblock {\"U}ber die {A}bgrenzung der {E}igenwerte einer {M}atrix.
\newblock {\em Bulletin de l'Acad\'emie des Sciences de l'URSS. Classe des
  sciences math\'ematiques et naturelles}, (6):749--754, 1931.

\bibitem{Goerisch:stufenverfahren-eigenwerten}
F.~Goerisch.
\newblock Ein {S}tufenverfahren zur {B}erechnung von {E}igenwertschranken.
\newblock In {\em Numerical treatment of eigenvalue problems, {V}ol. 4
  ({O}berwolfach, 1986)}, volume~83 of {\em Internat. Schriftenreihe Numer.
  Math.}, pages 104--114. Birkh\"{a}user, Basel, 1987.

\bibitem{Golberg-Chen:mfs-survey}
M.~A. Golberg and C.~S. Chen.
\newblock The method of fundamental solutions for potential, {H}elmholtz and
  diffusion problems.
\newblock In {\em Boundary integral methods: numerical and mathematical
  aspects}, volume~1 of {\em Comput. Eng.}, pages 103--176. WIT Press/Comput.
  Mech. Publ., Boston, MA, 1999.

\bibitem{GomezSerrano:survey-cap-in-pde}
J.~G\'{o}mez-Serrano.
\newblock Computer-assisted proofs in {PDE}: a survey.
\newblock {\em SeMA J.}, 76(3):459--484, 2019.

\bibitem{Gopal-Trefethen:new-laplace-solver-pnas}
A.~Gopal and L.~N. Trefethen.
\newblock New {L}aplace and {H}elmholtz solvers.
\newblock {\em Proc. Natl. Acad. Sci. USA}, 116(21):10223--10225, 2019.

\bibitem{Gopal-Trefethen:laplace-solver-detailed}
A.~Gopal and L.~N. Trefethen.
\newblock Solving {L}aplace {P}roblems with {C}orner {S}ingularities via
  {R}ational {F}unctions.
\newblock {\em SIAM J. Numer. Anal.}, 57(5):2074--2094, 2019.

\bibitem{Gordon-Webb-Wolpert:isospectral-domains-surfaces-riemannian-orbifolds}
C.~Gordon, D.~Webb, and S.~Wolpert.
\newblock Isospectral plane domains and surfaces via {R}iemannian orbifolds.
\newblock {\em Invent. Math.}, 110(1):1--22, 1992.

\bibitem{Grieser-Maronna:isospectral-triangle}
D.~Grieser and S.~Maronna.
\newblock Hearing the shape of a triangle.
\newblock {\em Notices Amer. Math. Soc.}, 60(11):1440--1447, 2013.

\bibitem{Guillemin-Kazhdan:inverse-spectral-negatively-curved-2-manifolds}
V.~Guillemin and D.~Kazhdan.
\newblock Some inverse spectral results for negatively curved {$2$}-manifolds.
\newblock {\em Topology}, 19(3):301--312, 1980.

\bibitem{Henrot:book-spectral-geometry}
A.~Henrot.
\newblock {\em Extremum problems for eigenvalues of elliptic operators}.
\newblock Frontiers in Mathematics. Birkh\"{a}user Verlag, Basel, 2006.

\bibitem{Hezari-Lu-Rowlett:neumann-isospectral-trapezoids}
H.~Hezari, Z.~Lu, and J.~Rowlett.
\newblock The {N}eumann isospectral problem for trapezoids.
\newblock {\em Ann. Henri Poincar\'{e}}, 18(12):3759--3792, 2017.

\bibitem{Hezari-Zelditch:inverse-spectral-problem-Rn}
H.~Hezari and S.~Zelditch.
\newblock Inverse spectral problem for analytic {$(\Bbb Z/2\Bbb
  Z)^n$}-symmetric domains in {$\Bbb R^n$}.
\newblock {\em Geom. Funct. Anal.}, 20(1):160--191, 2010.

\bibitem{Hezari-Zelditch:inverse-spectral-problem-ellipses}
H.~Hezari and S.~Zelditch.
\newblock One can hear the shape of ellipses of small eccentricity.
\newblock {\em Arxiv preprint arXiv:1907.03882}, 2019.

\bibitem{Judge-Hillariet:spectral-simplicity}
L.~Hillairet and C.~Judge.
\newblock Spectral simplicity and asymptotic separation of variables.
\newblock {\em Communications in Mathematical Physics}, 302(2):291--344, 2011.

\bibitem{Johansson:Arb}
F.~Johansson.
\newblock Arb: efficient arbitrary-precision midpoint-radius interval
  arithmetic.
\newblock {\em IEEE Transactions on Computers}, 66:1281--1292, 2017.

\bibitem{Judge-Mondal:hot-spots-triangles}
C.~Judge and S.~Mondal.
\newblock Euclidean triangles have no hot spots.
\newblock {\em Annals of Mathematics}, 191(1):167--211, 2020.

\bibitem{Kac:can-one-hear-shape-drum}
M.~Kac.
\newblock Can one hear the shape of a drum?
\newblock {\em Amer. Math. Monthly}, 73(4, part II):1--23, 1966.

\bibitem{Kato:upper-lower-bounds-eigenvalues}
T.~Kato.
\newblock On the upper and lower bounds of eigenvalues.
\newblock {\em J. Phys. Soc. Japan}, 4:334--339, 1949.

\bibitem{Laugesen-Siudeja:maximizing-neumann-triangles}
R.~S. Laugesen and B.~A. Siudeja.
\newblock Maximizing {N}eumann fundamental tones of triangles.
\newblock {\em J. Math. Phys.}, 50(11):112903, 18, 2009.

\bibitem{Laugesen-Siudeja:triangles-survey}
R.~S. Laugesen and B.~A. Siudeja.
\newblock Triangles and other special domains.
\newblock In {\em Shape optimization and spectral theory}, pages 149--200. De
  Gruyter Open, Warsaw, 2017.

\bibitem{Lehmann:optimale-eigenwerte}
N.~J. Lehmann.
\newblock Optimale {E}igenwerteinschliessungen.
\newblock {\em Numer. Math.}, 5:246--272, 1963.

\bibitem{Liu:framework-verified-eigenvalues}
X.~Liu.
\newblock A framework of verified eigenvalue bounds for self-adjoint
  differential operators.
\newblock {\em Appl. Math. Comput.}, 267:341--355, 2015.

\bibitem{Liu-Oishi:verified-eigenvalues-laplacian-polygons}
X.~Liu and S.~Oishi.
\newblock Verified eigenvalue evaluation for the {L}aplacian over polygonal
  domains of arbitrary shape.
\newblock {\em SIAM J. Numer. Anal.}, 51(3):1634--1654, 2013.

\bibitem{Lu-Rowlett:spectral-amm}
Z.~Lu and J.~Rowlett.
\newblock The sound of symmetry.
\newblock {\em Amer. Math. Monthly}, 122(9):815--835, 2015.

\bibitem{Miranda:PoincareMiranda}
C.~Miranda.
\newblock Un'osservazione su un teorema di {B}rouwer.
\newblock {\em Boll. Un. Mat. Ital. (2)}, 3:5--7, 1940.

\bibitem{Moler-Payne:bounds-eigenvalues}
C.~B. Moler and L.~E. Payne.
\newblock Bounds for eigenvalues and eigenvectors of symmetric operators.
\newblock {\em SIAM J. Numer. Anal.}, 5:64--70, 1968.

\bibitem{Nakao-Plum-Watanabe:cap-for-pde-book}
M.~T. Nakao, M.~Plum, and Y.~Watanabe.
\newblock {\em Numerical Verification Methods and Computer-Assisted Proofs for
  Partial Differential Equations}.
\newblock Springer, 2020.

\bibitem{Osgood-Phillips-Sarnak:compact-isospectral-plane-domains}
B.~Osgood, R.~Phillips, and P.~Sarnak.
\newblock Compact isospectral sets of plane domains.
\newblock {\em Proc. Nat. Acad. Sci. U.S.A.}, 85(15):5359--5361, 1988.

\bibitem{Osgood-Phillips-Sarnak:compact-isospectral-surfaces}
B.~Osgood, R.~Phillips, and P.~Sarnak.
\newblock Compact isospectral sets of surfaces.
\newblock {\em J. Funct. Anal.}, 80(1):212--234, 1988.

\bibitem{Osgood-Phillips-Sarnak:compact-isospectral-plane-domains-annals}
B.~Osgood, R.~Phillips, and P.~Sarnak.
\newblock Moduli space, heights and isospectral sets of plane domains.
\newblock {\em Ann. of Math. (2)}, 129(2):293--362, 1989.

\bibitem{Ourmieres:flat-triangles}
T.~Ourmi{\`e}res-Bonafos.
\newblock Dirichlet eigenvalues of asymptotically flat triangles.
\newblock {\em Asymptotic Analysis}, 92(3-4):279--312, 2015.

\bibitem{Parlett:symmetric-eigenvalue-book}
B.~N. Parlett.
\newblock {\em The symmetric eigenvalue problem}, volume~20 of {\em Classics in
  Applied Mathematics}.
\newblock Society for Industrial and Applied Mathematics (SIAM), Philadelphia,
  PA, 1998.
\newblock Corrected reprint of the 1980 original.

\bibitem{Paternain-Salo-Uhlmann:spectral-rigidity-anosov-surfaces}
G.~P. Paternain, M.~Salo, and G.~Uhlmann.
\newblock Spectral rigidity and invariant distributions on {A}nosov surfaces.
\newblock {\em J. Differential Geom.}, 98(1):147--181, 2014.

\bibitem{Plum:eigenvalues-homotopy-method}
M.~Plum.
\newblock Eigenvalue inclusions for second-order ordinary differential
  operators by a numerical homotopy method.
\newblock {\em Z. Angew. Math. Phys.}, 41(2):205--226, 1990.

\bibitem{Read-Sneddon-Bode:series-method-mps}
W.~W. Read, G.~E. Sneddon, and L.~Bode.
\newblock A series method for the eigenvalues of the advection diffusion
  equation.
\newblock {\em ANZIAM J.}, 45((C)):C773--C786, 2003/04.

\bibitem{Rellich:darstellung-eigenwerte}
F.~Rellich.
\newblock Darstellung der {E}igenwerte von {$\Delta u+\lambda u=0$} durch ein
  {R}andintegral.
\newblock {\em Math. Z.}, 46:635--636, 1940.

\bibitem{Still:computable-bounds-eigenvalues}
G.~Still.
\newblock Computable bounds for eigenvalues and eigenfunctions of elliptic
  differential operators.
\newblock {\em Numer. Math.}, 54(2):201--223, 1988.

\bibitem{Szabo:isospectral-metrics}
Z.~I. Szab\'{o}.
\newblock Isospectral pairs of metrics on balls, spheres, and other manifolds
  with different local geometries.
\newblock {\em Ann. of Math. (2)}, 154(2):437--475, 2001.

\bibitem{Tucker:validated-numerics-book}
W.~Tucker.
\newblock {\em Validated numerics}.
\newblock Princeton University Press, Princeton, NJ, 2011.
\newblock A short introduction to rigorous computations.

\bibitem{Watanabe:plane-domains-spectrally-determined}
K.~Watanabe.
\newblock Plane domains which are spectrally determined. {II}.
\newblock {\em J. Inequal. Appl.}, 7(1):25--47, 2002.

\bibitem{Weinstein-Stenger:eigenvalues-book}
A.~Weinstein and W.~Stenger.
\newblock {\em Methods of intermediate problems for eigenvalues}.
\newblock Academic Press, New York-London, 1972.
\newblock Theory and ramifications, Mathematics in Science and Engineering,
  Vol. 89.

\bibitem{Zelditch:inverse-spectral-problem-analytic-annals}
S.~Zelditch.
\newblock Inverse spectral problem for analytic domains. {II}.
  {$\mathbb{Z}_2$}-symmetric domains.
\newblock {\em Ann. of Math. (2)}, 170(1):205--269, 2009.

\end{thebibliography}

\end{document}